\documentclass[review,onefignum,onetabnum]{siamart171218_nolineno}

\usepackage{lipsum}
\usepackage{amsfonts}
\usepackage{graphicx}
\usepackage{epstopdf}
\usepackage{algorithm}
\usepackage{algpseudocode}
\usepackage{algorithmicx}
\usepackage{hyperref}
\usepackage{cleveref}

\newcounter{algsubstate}

\newenvironment{algsubstates}
{\setcounter{algsubstate}{0}%
	\renewcommand{\State}{%
		\refstepcounter{algsubstate}%
		\Statex {\footnotesize\alph{algsubstate}:}\space}}
{}
\usepackage{moreverb}
\usepackage{amsmath, amsfonts, amssymb,mathrsfs}
\usepackage{algorithm}
\usepackage{algpseudocode}
\usepackage{bm}
\usepackage{listings}
\usepackage{multirow}

\ifpdf%
  \DeclareGraphicsExtensions{.eps,.pdf,.png,.jpg}
\else
  \DeclareGraphicsExtensions{.eps}
\fi
\usepackage{amsopn}

\usepackage{booktabs}

\newcommand{\TheTitle}{%
  An Adaptive Multigrid Method Based on Path Cover
}

\newcommand{\TheShortTitle}{%
  An Adaptive AMG Based on Path Cover
}

\newcommand{\TheName}{%
  Junyuan Lin
}

\author{Xiaozhe Hu \thanks{Department of Mathematics, Tufts
		University, Medford, MA 02155 (\hbox{xiaozhe.hu@tufts.edu},
		\hbox{junyuan.lin@tufts.edu}).} 
	\and 
	Junyuan Lin 
	\footnotemark[2] 
	\and
	Ludmil T. Zikatanov \thanks{Department of Mathematics, The Pennsylvania
		State University, University Park, PA 16802 (\hbox{ludmil@psu.edu}).}
}
\title{{\TheTitle}\thanks{Submitted to the editors on June 15, 2018.
		\funding{The work of Hu and Lin is partially supported by NSF grant DMS-1620063. The work of Zikatanov is supported in part by NSF DMS-1720114.}}}
\headers{\TheShortTitle}{X.~Hu, J.~Lin, and L.~T.~Zikatanov}
\ifpdf%
\hypersetup{%
	pdftitle={\TheTitle},
	pdfauthor={\TheName}
}
\fi

\begin{document}
\baselineskip=0.94\normalbaselineskip
\maketitle


\begin{abstract}
We propose a path cover adaptive algebraic multigrid (PC-$\alpha$AMG) method for solving linear systems of weighted graph Laplacians and
can also be applied to discretized second order elliptic partial
differential equations. The PC-$\alpha$AMG is based on unsmoothed
aggregation AMG (UA-AMG). To preserve the structure of smooth error down to the coarse levels, we approximate the level sets of the smooth
error by first forming vertex-disjoint path cover with paths
following the level sets. The aggregations are then formed by matching along the paths in the path cover.  In such manner, we are able
to build a multilevel structure at a low computational cost.  The proposed PC-$\alpha$AMG provides a mechanism
to efficiently re-build the multilevel hierarchy during the iterations
and leads to a fast nonlinear multilevel algorithm.
Traditionally, UA-AMG requires more sophisticated cycling techniques,
such as AMLI-cycle or K-cycle, but as our numerical results show, the PC-$\alpha$AMG proposed here leads to nearly
optimal standard V-cycle algorithm for solving linear systems with
weighted graph Laplacians. Numerical experiments for some real world
graph problems also demonstrate PC-$\alpha$AMG's effectiveness and robustness,
especially for ill-conditioned graphs.
	\end{abstract}



\section{Introduction}\label{sec:intro}
As weighted graphs frequently being employed as the data representations to describe a rich spectrum of application fields, including social, sensor, energy, and neuronal networks~\cite{borgatti2009social,yang2010sensor,holmgren2006power,bullmore2009brain}, the associated graph Laplacians naturally arise in large-scale computations of various application domains. In recent works, solving weighted graph Laplacians has been applied to solve ranking and user recommendation problems~\cite{Hirani.A;Kalyanaraman.K;Watts.S2011a,Jiang.X;Lim.L;Yao.Y;Ye.Y2011a,HodgeRank}.  In~\cite{Cao.M;Pietras.C;Feng.X;Doroschak.K;Schaffner.T;Park.J;Zhang.H;Cowen.L;Hescott.B2014a,Cao.M;Zhang.H;Park.J;Daniels.N;Crovella.M;Cowen.L;Hescott.B2013a,DSD2018}, similarities of proteins are calculated by solving the graph Laplacian associated with the protein interaction network,  which are further used in clustering and labeling proteins' by their functions. The marriage of graph Laplacian and heated computer science topics such as Convolutional Neural Networks and tensor decomposition has also been a trend~\cite{he2006tensor,Bronstein.M;Bruna.J;LeCun.Y;Szlam.A;Vandergheynst.P2017a,Kipf.T;Welling.M2016a,Henaff.M;Bruna.J;LeCun.Y2015a,Bruna.J;Zaremba.W;Szlam.A;LeCun.Y2013a,shaham2018spectralnet,li2018adaptive}. Researchers advanced algorithms that adapt graph Laplacian properties to improve tasks including image reconstruction, clustering image data-sets and classification~\cite{agarwal2006higher,narita2012tensor,GLTD,steerableGL}. To efficiently solve graph Laplacians, algebraic multigrid (AMG) is often applied.  The standard AMG method was proposed to solve partial differential equations (PDEs) and involves mainly two parts: smoothing out the high-frequency errors on the fine levels and eliminating the low-frequency errors on the coarse grids~\cite{xu2017algebraic,Vassilevski.P2008,Ruge.J;Stuben.K1987,AMG_1984,ruge1984efficient,ruge1983algebraic}. AMG is proven to be one of the most successful iterative methods in practical applications and many AMG methods have been developed to solving graph Laplacians, such as combinatorial multigrid~\cite{koutis2011combinatorial}, Lean AMG~\cite{livne2012lean}, Algebraic multilevel preconditioner based on matchings/aggregations~\cite{Kim.H;Xu.J;Zikatanov.L2003,notay2010aggregation,Brannick.J;Chen.Y;Kraus.J;Zikatanov.L.2013a,d2013adaptive}.


Traditional AMG methods usually build the multilevel structures beforehand and there is no interplay between the remaining errors and coarsening schemes. This results in many computational inefficiencies. For example, the cycling scheme applies the same multilevel hierarchy during the solve phase even when convergence rate is slow, which leads to computational waste.  To better reduce the computational redundancy, adaptivity becomes essential.
%
Substantial efforts have been made to incorporate adaptivity in iterative methods. 
Date back to 1984, the original adaptive AMG algorithm~\cite{AMG_1984} was proposed, which laid foundation for the development of self-learning and bootstrap AMG. The bootstrap AMG approach was further developed in~\cite{manteuffel2010operator, brandt2011bootstrap} which starts with several test vectors while the adaptive approaches typically start with only one test vector. Directed by the theory of smoothed aggregation AMG (SA-AMG) developed in~\cite{vanvek1996algebraic,van2001convergence}, Brezina et al. introduced adaptive SA-AMG that determines interpolation operators based on information from the system itself rather than on explicit knowledge of the near-kernel space~\cite{brezina2005adaptive}. In the following year, the same authors further proposed an operator-induced interpolation approach that automatically represent smooth components~\cite{brezina2006adaptive}. In both of the literatures, the basic idea is to solve the homogeneous problem ($\bm{A}\bm{x}=\bm{0}$) to expose the slow-to-converge errors, coarsening based on algebraically smooth errors, and then applying back to improve the AMG process. MacLachlan et al. further developed a two-level, reduction-based nonlinear adaptive AMG and achieved local convergence and possibly global convergence in special cases~\cite{maclachlan2006adaptive}. The development of the UA-AMG method is fairly recent.  In~\cite{d2013adaptive}, D'Ambra et al. applied coarsening scheme that uses compatible weighted matching, which avoids the reliance on the characteristics of connection strength in defining both the coarse space and the interpolation operators. Through previous endeavors, authors enlightened different ideas on how to build smooth errors into the ranges of adaptively constructed interpolation operators, assuming coarsening is done. This motivates us to the idea of using the smooth error for not only building interpolation operators but also determining coarsening scheme.

In this paper, an adaptive AMG method based on path cover is presented. Specifically, we choose to use UA-AMG~\cite{Blaheta.R1988} for its low computational complexity, which is favorable for parallel computing. We setup the multilevel structure by first applying path cover~\cite{ore1961arc} to approximate the level sets of the smooth error, then coarsening along the paths so that the structure of smooth error is well represented on the coarse levels, and building the smooth error into the range of interpolation operator to effectively eliminate it.  The setup phase is reactivated when the convergence rate becomes slow, which makes the scheme adaptive to the slow-to-converge errors.  Instead of using AMLI- or K-cycles for standard UA-AMG~\cite{Brannick.J;Chen.Y;Kraus.J;Zikatanov.L.2013a,Hu.X;Vassilevski.P;Xu.J2012a,Vassilevski.P2008}, we simply use V-cycle and achieve nearly uniform convergence for model problems. The design of combining UA-AMG and V-cycle has advantages in its simplicity and efficiency, therefore is useful for parallel computing, especially on solving matrices with large condition numbers. While the existed adaptive methods need to first solve the homogeneous problem, $\bm{A}\bm{x}=\bm{0}$, to determine the near-null-space components and use this knowledge to solve $\bm{A}\bm{x}=\bm{b}$, our adaptive method solves $\bm{A}\bm{x}=\bm{b}$ directly, for a general right-hand-side $\bm{b}$.

The paper is arranged as follows. In \cref{sec:pre}, basic subroutines, such as a modified version of path cover finding algorithm, are presented.  Our main adaptive AMG algorithm is discussed in \cref{sec:alg}.  Numerical results are demonstrated and analyzed in \cref{sec:num}. Finally, in \cref{sec:conc}, we conclude the main contributions of this paper and enumerate possible future directions.

\section{Preliminaries}\label{sec:pre}
In this section, we review basic aggregating methods and more importantly, the path cover approximation algorithm~\cite{moran1990approximation}. We use these components in our main algorithm presented in \cref{sec:alg}. 
%
%
%
%
\subsection{Basic Algorithms for Aggregation}
There are various off-the-shelf aggregating methods that are frequently applied to graphs, such as maximal independent set~\cite{MIS1986}, heavy edge coarsening algorithm~\cite{urschel2014cascadic} and maximal weighted matching (MWM)~\cite{galil1986ev}. In this paper, we choose to compare the performance of our matching-like aggregating method (\cref{alg:matching}) with the most closely related MWM. 
Similar comparison results could be observed between \cref{alg:matching} and other graph matching methods as well. Here, we briefly recall MWM, since in the implementations, we adopt this aggregating scheme on weighted adjacency matrix $\bm{A}$ of each model graph as a default and set the performance benchmark of regular AMG method.

MWM algorithm forms matchings by visiting edges in the graph, from heaviest to lightest, and matches two endpoints of the edge if they are unaggregated. Commonly, there might be isolated nodes after applying MWM.  We add those isolated nodes to existing matchings in order to keep the number of aggregates low, and at the same time we set three as an upper bound to the diameter of each aggregate. 

\subsection{Constructing a Vertex Disjoint Path Cover}
Consider a graph $G = (V,E,\omega)$ with positive weights but no self-loops or parallel edges, where $V$ is the vertex set, $E$ is the edge set, and $\omega > 0$ represents the weight set.   A \textit{path cover}, $S$, of $G$ is a set of vertex-disjoint paths that covers all the vertices of $G$~\cite{ore1961arc}. Now define
\begin{equation*}
\omega^{\ast}(G):=\max_{S}\sum_{e\in S}\omega(e),
\end{equation*}
which is the maximum possible weight of a path cover $S$ of $G$. In~\cite{moran1990approximation}, Moran et al. pointed out that to find the exact maximal weighted path cover of a graph is an NP-complete problem, but on the other hand, they showed that one can find a $\frac{1}{2}$-approximated path cover, $\widetilde{S}$ of $G$, in $\mathcal{O}(|E|\cdot\log|E|)$ time, such that
\begin{equation*}
\omega_{\widetilde{S}}(G)=\sum_{e\in\widetilde{S}}\omega(e)\geqslant\frac{1}{2} \omega^{\ast}(G).
\end{equation*}
The $\frac{1}{2}$-Approximation path cover finding algorithm presented in \cref{alg:PC} is a slightly modified version of the algorithm from~\cite{moran1990approximation}.

\begin{figure}[h!]
	\begin{center}
		\includegraphics[width=0.5\textwidth]{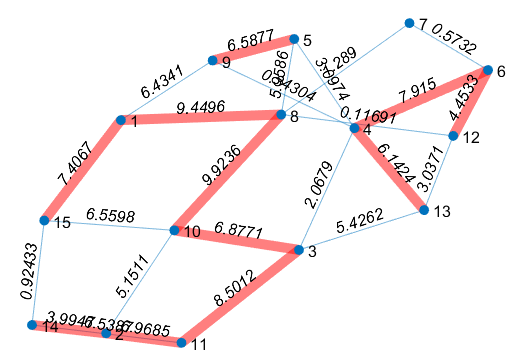}	
	\end{center}
	\caption{Path cover found by \cref{alg:PC} on a random graph with weights}
	\label{fig:path-cover}
\end{figure}

\begin{algorithm}
	\caption{Path Cover}\label{alg:PC}
	\begin{algorithmic}[1]
		\Procedure{[ $\texttt{cover}$ ] = \texttt{PathCover}($\bm{A} $)}{}
		\State Input:
		$\bm{A} $ - adjacency matrix of an undirected positive weighted graph $G=(V,E, \omega)$
		\State Output:
		$\texttt{cover}$ - a path cover of graph $G$ 
		\Statex
		\State $\texttt{sorted\_edges} \gets$ \textbf{Sort} the edges in descending order  based on weights.
		
		\For{$e(u,v)\in \texttt{sorted\_edges}$ }
		\If {both $u$ and $v$ are endpoints}
		\If {neither $u$ nor $v$ is in any paths}
		\State 
		Add $\{u,v\}$ as a new path in $\texttt{cover}$
		\ElsIf{$u$ is in a path, but $v$ is not}
		\State
		Append $\{v\}$ to $\texttt{cover}\{\text{path that contains }u\}$
		\ElsIf{$v$ is in a path, but $u$ is not}
		\State
		Append $\{u\}$ to $\texttt{cover}\{\text{path that contains }v\}$
		\Else 
		\Comment ($u$ and $v$ are in different nonzero-length paths)
		\State
		Merge two paths
		\EndIf
		\EndIf
		\EndFor
		
		\EndProcedure
		
	\end{algorithmic}
\end{algorithm}
\cref{alg:PC} is a greedy algorithm that checks each edge, from the heaviest to the lightest, and incorporates edges as long as the set still contains only paths. 
The complexity of \cref{alg:PC} is $\mathcal{O}(|V|\log |V|)$ when the graph is sparse, $\mathcal{O}(|E|)=\mathcal{O}(|V|)$. \cref{fig:path-cover} illustrates a resulting path cover of a random graph, with weights displayed on its edges. Notice that there might be isolated points after finding the path cover. Therefore, later we include the isolated points using \cref{alg:matching}, so that all the points on the original graph are represented on coarse levels.

\section{Adaptive Algorithm}\label{sec:alg}
The main idea and details of our PC-$\alpha$AMG are presented in this section.  In \cref{sec:alg_pre}, the main idea of the
algorithm is presented and demonstrated with an intuitive
example. \cref{sec:alg_approx} and \cref{sec:alg_PC}
comprise essential subroutines of approximating smooth error using
previously built multilevel hierarchies
(\cref{alg:ApproxSE}), building aggregations and prolongation
operator from path cover (\cref{alg:matching}), and setting
up multilevel structure using path cover and approximated error
(\cref{alg:pcsetup}). Finally, the full PC-$\alpha$AMG algorithm is presented in
\cref{sec:alg_main}.
\subsection{Basic Ideas and Rationale of the Adaptive Algorithm}\label{sec:alg_pre}
As previously mentioned, AMG reduces error during two procedures:
smoothing out the high-frequency errors and eliminating low-frequency
errors that are restricted to the coarse grids. Our adaptive algorithm
contributes to the latter aspect by preserving the smooth error well
onto the coarse levels, so that the adaptively designed multilevel
hierarchy can eliminate the current smooth error efficiently on
the coarsest level. Specifically, we utilize path cover finding
algorithm to capture the level sets of the smooth error on the finest
level, aggregate/match along these paths to preserve the smooth error
on the coarser levels, and reconstruct AMG hierarchy to aim at this
specific remaining smooth error. In this manner, we can eliminate the
dominating smooth errors which cause slow convergence one by one,
until desired accuracy is met.  
  In our opinion, it could be beneficial
to have multilevel hierarchies which approximate the errors well when
$\bm{b}\neq 0$.  One possible approach, albeit beyond the scope of our
considerations here, is to use adaptive aggregations based on a
posteriori error estimates on graphs as proposed in~\cite{Wenfang}.
Unlike other adaptive AMG methods, PC-$\alpha$AMG algorithm proposed here integrates the setup and solve
phase together by identifying the smooth
errors while solving the linear system $\bm{A}\bm{x}=\bm{b}$.

\begin{figure}[h!]
	\begin{center}
		\includegraphics[width=0.3\textwidth]{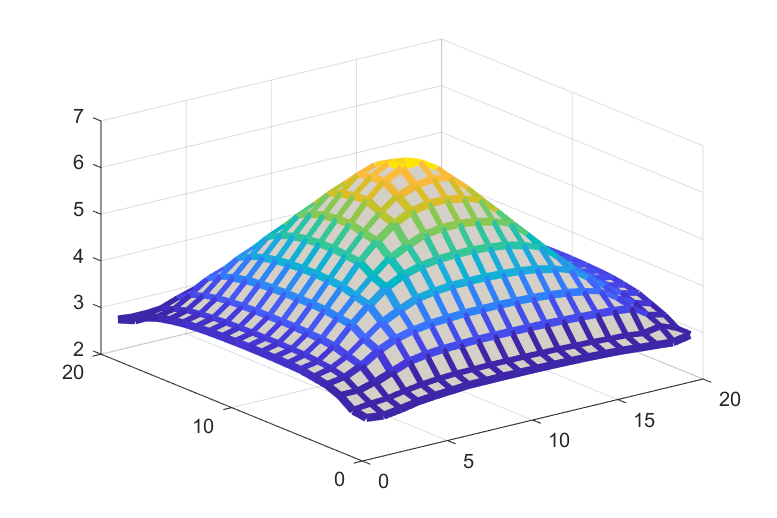} \ \ 
		\includegraphics[width=0.3\textwidth]{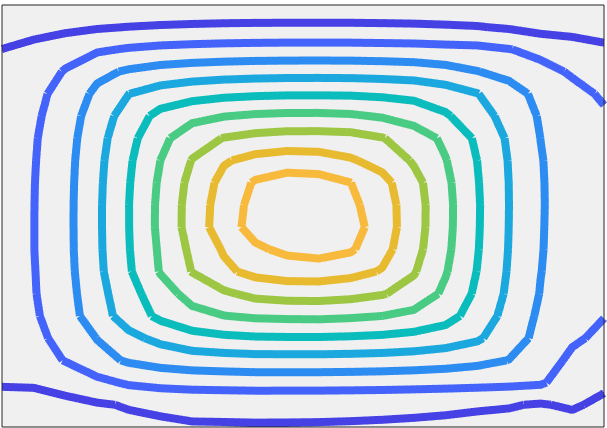} \\	
		\includegraphics[width=0.3\textwidth]{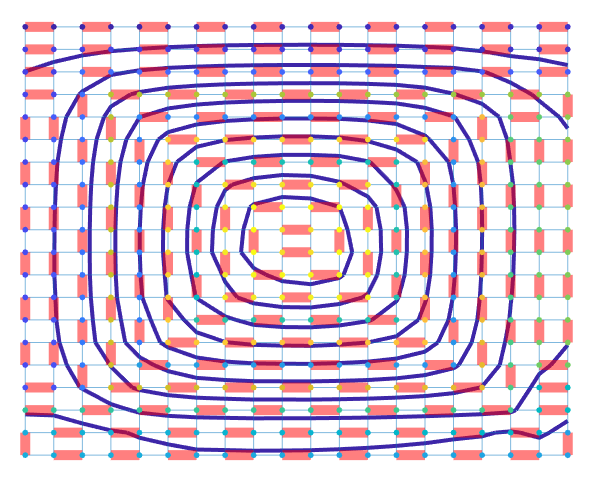}  \ \
		\includegraphics[width=0.3\textwidth]{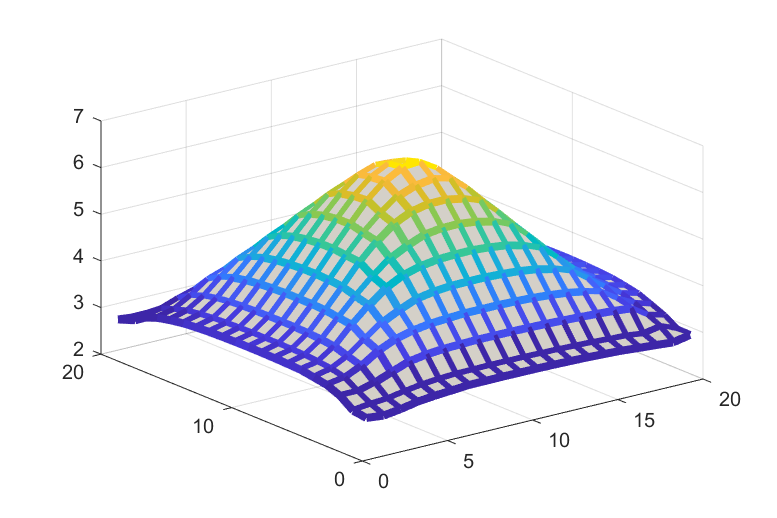}
	\end{center}
	\caption{Upper Left: smooth error; Upper Right: level sets of smooth error; Lower Left: matchings on level sets; Lower Right: error after restricting and prolongating}
	\label{fig:intuition}
\end{figure}

To test if the aggregations/matchings along the level sets of smooth error can successfully preserve the smooth error onto coarse levels, we take a manufactured smooth error (\cref{fig:intuition} upper left), manually build matchings (\cref{fig:intuition} lower left) on its level sets (\cref{fig:intuition} upper right), and use the matchings to build the prolongation operator to restrict and prolongate the error back to the original level (\cref{fig:intuition} lower right). The norm of difference between the two smooth errors is within $10^{-10}$. This reassures us that the aggregating scheme based on level set is efficient in capturing the smooth error.  In our algorithm, we use path cover to approximate level sets, where each path in the path cover represents one level set.  In this way, finding level sets can be done purely algebraically using only the matrices.
In \cref{sec:alg_approx}, we include procedures that can obtain the estimated smooth error when solving slows down, which is symbolized by the manufactured error (\cref{fig:intuition} upper left). We present subroutines in \cref{sec:alg_PC} to automatically aggregate along each path in path cover that approximates the level sets of the approximated smooth error. It is verifiable that we can achieve similar effect as the manufactured aggregates (\cref{fig:intuition} lower left) have in preserving the smooth error structure.

\subsection{Approximating the Smooth Error}\label{sec:alg_approx}
Since all the subroutines that aim to preserve and ultimately eliminate the smooth error on the coarse levels rely on the fact that the smooth error itself can be well approximated on the fine level, one essential step of our algorithm is to approximate the smooth error accurately when the convergence rate becomes slow. The fact is, during the iterations, the approximated solution $\bm{x}_k$ is available at hands, but not straight-forwardly the error $\bm{e}_k$. 
Based on the well-known error equation $\bm{A}\bm{e}_k = \bm{r}_k = \bm{b} - \bm{A} \bm{x}_k$, we try to approximately solve the error equation by applying several steps of multigrid preconditioned conjugate gradient (PCG) method and detailed approach is presented in \cref{alg:ApproxSE}.  For the first several times of re-setups,  we use several W-cycles to get a good approximation (this corresponds to the case ``$\texttt{Re} \leqslant \texttt{numW}$"). In practice, we minimized the usage of W-cycle (see \cref{sec:num} for details).  For later re-setups, we simply use V-cycles based on the existing multilevel structures.  More precisely,  in Step 10 of \cref{alg:ApproxSE}, 
 we use a symmetric composite preconditioner proposed in~\cite{d2013adaptive}:
\begin{align*}
	\bm{e}\gets\prod_{j=0}^{2\cdot\texttt{Re}+1}(\bm{I}-\bm{B}_j\bm{A})\bm{e},
\end{align*}
where \texttt{Re} is the number of re-setups, and $\bm{B}_{\texttt{Re}+j} = \bm{B}_{\texttt{Re}+1-j}, j = 1,\dots \texttt{Re}+1$. Each $\bm{B}_j$ corresponds to the preconditioning effect of multilevel hierarchy built from the $j$-th re-setup, i.e. $\texttt{P\_hist}\{j\}$. In this manner, we recycle all the hierarchies ever built and reduce the computational waste.
\begin{algorithm}[h!]
	\caption{Approximate Smooth Error}\label{alg:ApproxSE}
	\begin{algorithmic}[1]
		\Procedure{[ $\bm{e}$ ] = \texttt{ApproximateSmoothError}($\bm{A},\{\bm{A}_{\ell}\}^L_{\ell=2},\bm{r}, \texttt{P\_hist}, \texttt{Re}, \bm{e}$)}{}
		\State
		Inputs: 
		\begin{algsubstates}
			\State $\bm{r}$ - residual vector
			\State \texttt{P\_hist} - all the hierarchies ever created
			\State \texttt{Re} - re-setup count
		\end{algsubstates}
		\Statex
		\State
		Choose parameters: 
		\begin{algsubstates}
			\State \texttt{numW}
			\Comment Number of smooth errors (re-setups) using W-cycle
			\State \texttt{iterW}
			\Comment Number of iterations of W-cycle PCG
			\State \texttt{iterV}
			\Comment Number of iterations of V-cycle PCG
		\end{algsubstates}
		\If {$\texttt{Re} \leqslant \texttt{numW}$}
			\For {$i \gets 1 : \texttt{iterW}$}
				\State $\bm{e}\gets\texttt{W\_cycle\_PCG}(\bm{A},\bm{r},\bm{e},\texttt{P\_hist\{\texttt{Re}\}},\{\bm{A}_{\ell}\}^L_{\ell=2})$
			\EndFor
		\Else
			\For {$i = 1 : \texttt{iterV}$}
			\State $\bm{e}\gets\texttt{V\_cycle\_PCG}(\bm{A},\bm{r},\bm{e},\texttt{P\_hist\{\texttt{1:Re}\}},\{\bm{A}_{\ell}\}^L_{\ell=2})$
			\EndFor
		\EndIf
		\State Normalize $\bm{e}$
		\State Make $\bm{e}$ orthogonal to $\bm{1}$   \Comment{This is because we are solving Graph Laplacian}
		\EndProcedure
		
	\end{algorithmic}
\end{algorithm}

\subsection{Coarsening based on Path Cover}\label{sec:alg_PC}

\begin{algorithm}[h]
	\caption{Path Cover Aggregation}\label{alg:matching}
	\begin{algorithmic}[1]
		\Procedure{[ $\bm{P}$ ] = \texttt{PathCoverAggregate}(\texttt{cover},$\widetilde{\bm{A}},\ell, \bm{e}$)}{}
		\State Input:
		\begin{algsubstates}
			\State $\widetilde{\bm{A}}$ - graph Laplacian (possibly reweighted)
			\State $\ell$ - level count
			\State $\bm{e}$ - smooth error
		\end{algsubstates}
		\State Output:
		$\bm{P}$ - prolongation operator
		\Statex
		\State
		$n\gets$ number of nodes
		\For{each path in \texttt{cover}}
		\State
		Match each vertex with one of its neighbor
		\EndFor
		\State
		$\texttt{numAgg}\leftarrow$ number of paths in \texttt{cover}
		\For{each unaggregated vertex $u$}  
		\State
		$v\gets$ $u$'s neighbor with the maximal edge weight based on $\widetilde{\bm{A}}$
		\If{$v$ is also unaggregated}
		\Comment Case 1
		\State
		Match $u$ and $v$ together
		\ElsIf {\texttt{size}(aggregation of $v$) $=2$}
		\Comment Case 2
		\State
		Aggregate $u$ into $v$'s aggregation
		\Else
		\Comment Case 3, \texttt{size}(aggregation of $v$) $=3$
		\State
		Match $u$ and $v$ and remove $v$ from previous aggregation of $v$
		\State
		$\texttt{numAgg} \leftarrow \texttt{numAgg} + 1$
		\EndIf 
		\EndFor
		\If {$\ell=1$}
		\State
		Build prolongation operator $\bm{P}\in \mathbb{R}^{n\times\texttt{numAgg}}$ using smooth error $\bm{e}$
		\Else
		\State
		Build boolean prolongation operator $\bm{P}\in \mathbb{R}^{n\times\texttt{numAgg}}$ 
		\EndIf
		\EndProcedure
		
	\end{algorithmic}
\end{algorithm}

\begin{algorithm}[h!]
	\caption{Path Cover AMG Setup}\label{alg:pcsetup}
	\begin{algorithmic}[1]
		\Procedure{[ $\{\bm{A}_{\ell}\}^L_{\ell=2},\{\bm{P}_{\ell}\}^L_{\ell=1}$ ] = \texttt{PathCoverAMG\_setup}($\bm{A},\bm{e}$)}{}
		
		\State
		$\widetilde{\bm{A}}_{ij} \gets  1/\|\bm{e}_i-\bm{e}_j\|$, if edge $e(i,j)$ is in the graph induced by $\bm{A}^2$
		\State
		$\texttt{cover} \gets$ \texttt{PathCover}($\widetilde{\bm{A}}$)
		\State
		$\bm{A}_1\gets \bm{A}$
		\For {$\ell \gets 1:L$}  
		\State
		$\bm{P}_{\ell} \gets \texttt{PathCoverAggregate}(\texttt{cover},\widetilde{\bm{A}},\ell, \bm{e})$
		\State
		$\widetilde{\bm{A}}\gets\bm{P}_{\ell}^T\widetilde{\bm{A}}\bm{P}_{\ell}$
		
		\State
		Shorten each path in \texttt{cover} by merging the vertices in one aggregate into a vertex
		\State
		$\bm{A}_{\ell+1}\gets \bm{P}_{\ell}^T \bm{A}_{\ell} \bm{P}_{\ell}$ 
		\EndFor 
		\EndProcedure
		
	\end{algorithmic}
\end{algorithm}

\begin{figure}[h!]
	\begin{center}
		\includegraphics[width=0.3\textwidth]{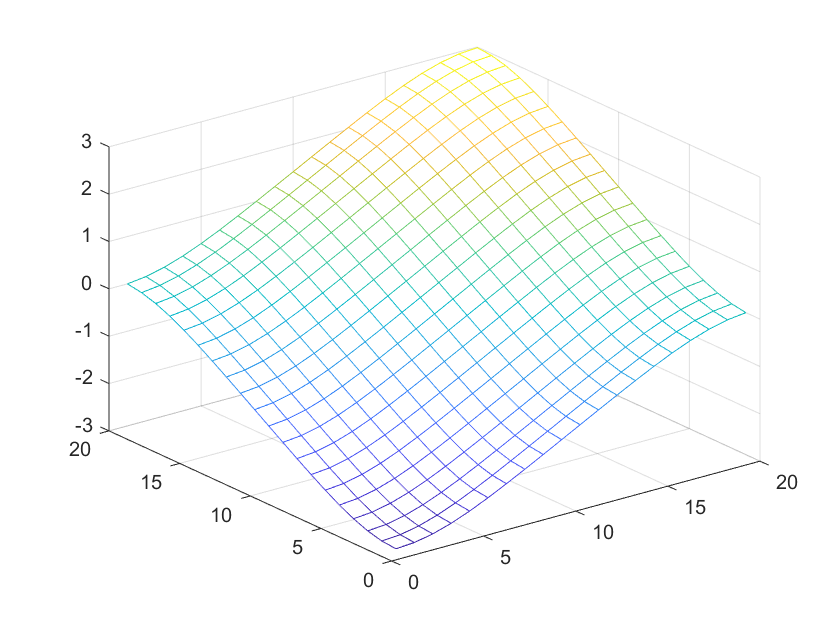}	
		\includegraphics[width=0.3\textwidth]{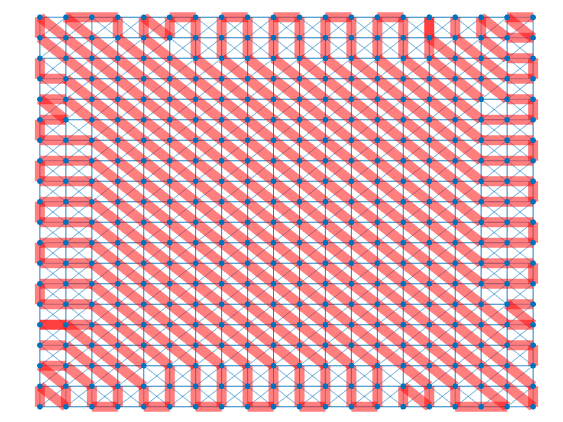}
		\includegraphics[width=0.3\textwidth]{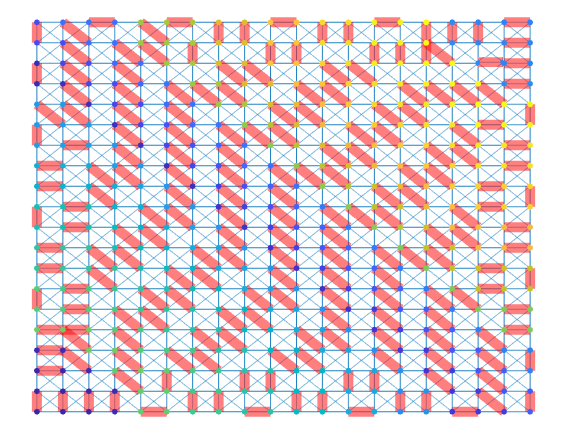}
	\end{center}
	\caption{Left: smooth error; Middle: path cover; Right: matching/aggregates on path cover}
	\label{fig:path-cover-matching}
\end{figure}


We use the information of the (approximate) smooth error in forming aggregations. In \cref{fig:path-cover-matching}, we illustrate the process of approximating the level sets of a smooth error (\cref{fig:path-cover-matching}, left) using path cover finding algorithm (\cref{alg:PC}) and aggregating along the path cover (\cref{alg:matching}).  We reassign the weights in \cref{alg:pcsetup}, Step 2, i.e. $\widetilde{\bm{A}}_{ij} =1/\|\bm{e}_i-\bm{e}_j\|$, so that edges between nodes of similar values (in another word, those nodes should be on the same level set) have heavy weights and are more likely to be chosen to form the path cover.
However, the edges in the original graph might not align with the heavy edges (e.g., the diagonal lines in the leftmost subplot in \cref{fig:path-cover-matching}). We therefore use the sparsity pattern of $\bm{A}^2$ to include more edges while keeping the sparsity of the graph relatively low (For the example shown in \cref{fig:path-cover-matching}, sparsity pattern of $\bm{A}^2$ includes the diagonal edges which are exactly on the level sets of the smooth error). After reassigning weights and changing the sparsity pattern, \cref{alg:PC} finds the path cover.  As shown in the middle subplot of \cref{fig:path-cover-matching} as an example, the path cover approximate level sets of the smooth error as we expected. Note that, since the smooth error are not perfectly linear near the boundaries, the level sets are not strictly diagonals.

\begin{figure}[h!]
	\begin{center}
		\includegraphics[width=0.9\textwidth]{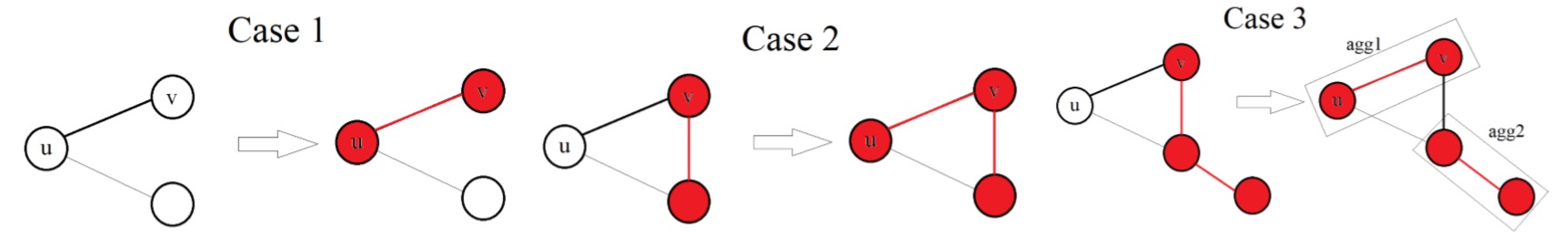}	
	\end{center}
	\caption{Aggregating the isolated points}
	\label{fig:Aggregating}
\end{figure}
As previously demonstrated, there might be isolated points resulting
from creating path cover and possibly other sources. 
\cref{alg:matching} handles the isolated points altogether.
\cref{fig:Aggregating} corresponds to Steps 9-19 in
\cref{alg:matching} and shows all three possible cases.  As we can
see in the right subplot of \cref{fig:path-cover-matching}, after rearrangement, there is no isolated vertex and the
resulting aggregations are all of size $2$ or $3$.


Here, \cref{fig:path-cover-matching} only explicitly shows the process of coarsening on the finest level. If the smooth error is well preserved onto the coarse level, its level sets should not change much on the coarse levels.  Thus, we simply match the nodes on each path in the fine level path cover to generate coarse-level path cover and recurse this process for all coarse levels. 

Also notice that we use the smooth error again in building prolongation operator in \cref{alg:matching}, namely, the prolongation operator $\bm{P}_{\ell}$ is defined as
$$
(\bm{P}_{\ell})_{ki} = 
\begin{cases}
\bm{e}_k & \text{if} \ k \in V^i_{\ell}, \ell=1, \\
1 & \text{if} \ k \in V^i_{\ell}, \ell>1, \\
0  & \text{otherwise},
\end{cases}
$$
where $\ell$ is the level count and $\bm{e}_k$ is the $k$th entry of the (approximate) smooth error $\bm{e}$. Obviously, we have~$\bm{e}\in \text{Range}(\bm{P}_1)$.


\subsection{Main Algorithm}\label{sec:alg_main}
\begin{figure}[h!]
	\begin{center}
		\includegraphics[width=0.7\textwidth]{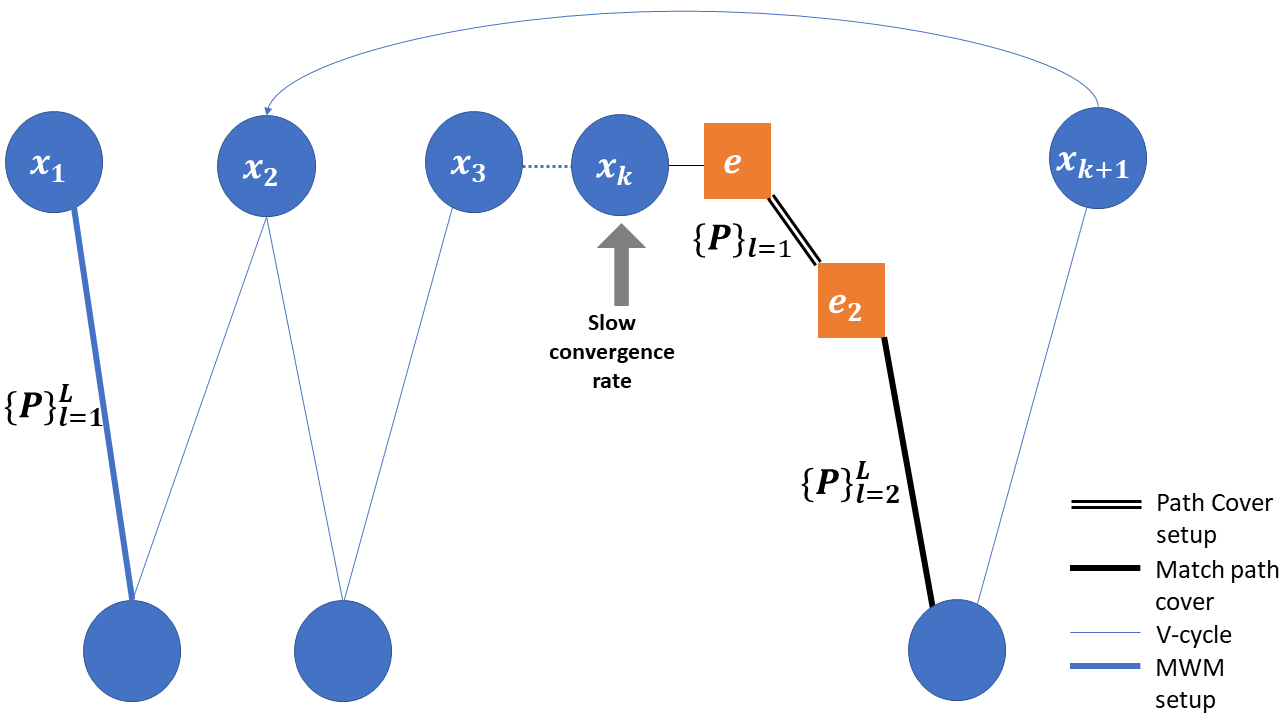}
		\caption{Procedure of PC-$\alpha$AMG}\label{fig:adpAMG}
	\end{center}
\end{figure}

The overall PC-$\alpha$AMG algorithms are present in \cref{alg:adpAMG_b} and \cref{alg:adpAMG_0}.  \cref{alg:adpAMG_b} is for the general case $\bm{A}\bm{x}=\bm{b}$. 
\cref{alg:adpAMG_0} is a special case for solving $\bm{A}\bm{x}=\bm{0}$. One reason we include the homogeneous case is that we have the knowledge that $\bm{x}=\bm{0}$ is the true solution in this case, so that the process of attaining smooth error is much simpler.  Another reason is that following the basic idea of adaptive AMG methods, \cref{alg:adpAMG_0} can also be used as standalone setup phase and build several multilevel structures that are effective for eliminating smooth errors.
The two versions of algorithms include all of the aforementioned components, yet the main procedure of both versions is rather straight-forward. We use MWM to build the initial multilevel structure and use it as the preconditioner to solve the model problem.  For simple graph Laplacian problems, MWM might already give good performance, then there is no need to use the adaptive procedure.  Our PC-$\alpha$AMG aims at difficult problems which have
error components that are slow-to-converge for MWM, the current multilevel hierarchy is not desirable anymore. We then re-setup using aggregations based on path cover of the current smooth error.  This is done whenever the convergence rate becomes higher than a certain pre-described threshold.  Such adaptive approach efficiently eliminates slow-to-converge smooth errors. We use \cref{fig:adpAMG} to demonstrate our PC-$\alpha$AMG.

\begin{algorithm}[h!]
	\caption{PC-$\alpha$AMG (for general $\bm{b}$)}\label{alg:adpAMG_b}
	\begin{algorithmic}[1]
		\Procedure{[ $\bm{x}$ ] = \texttt{AdaptiveAMG}($\bm{A},\bm{b},\bm{x},\texttt{tol},\texttt{max\_iter}, \texttt{threshold}$)}{}
		\State
		$\{\bm{A}_{\ell}\}^L_{\ell=2},\{\bm{P}_{\ell}\}^L_{\ell=1}\gets$\texttt{MWM\_setup}($\bm{A}$)
		\State
		$k \gets 1$
		\Comment Initialization for iteration number
		\State
		\texttt{Re} $\gets 1$
		\Comment Initialization for re-setup counts
		\State
		\texttt{P\_hist}$\{\texttt{Re}\} \gets \{\bm{P}_{\ell}\}^L_{\ell=1}$
		\Comment \texttt{P\_hist} saves all hierarchies ever created
		\State
		$\bm{e} \gets \bm{0}$
		\Comment Initialize approximated smooth error at $\bm{0}$
		\State
		$\bm{x}_k \gets \bm{x}$
		\State
		$\bm{r}_k \gets \bm{b}-\bm{A} \bm{x}_{k}$
		\While {$k<$ \texttt{max\_iter} and $\|\bm{r}_k\|>\texttt{tol}$} 
		\State
		$\bm{x}_{k+1}\gets$ \texttt{V\_cycle}$(\bm{A},\bm{b},\bm{x}_k,\{\bm{P}_{\ell}\}^L_{\ell=1},\{\bm{A}_{\ell}\}^L_{\ell=2})$
		\State 
		Make $\bm{x}_k$ orthogonal to $\bm{1}$
		\Comment If we are solving graph Laplacian
		\State
		$\bm{r}_{k+1}\gets \bm{b}-\bm{A} \bm{x}_{k+1}$
		\State
		$\texttt{ConvR}\gets\frac{\|\bm{r}_{k+1}\|}{\|\bm{r}_{k}\|}$
		\Comment Compute convergence rate at $(k+1)$ step
		\If{average \texttt{ConvR} of iterations after the last re-setup $>\texttt{threshold}$}
		\State
		$\bm{e}\gets\texttt{ApproximateSmoothError}(\bm{A},\{\bm{A}_{\ell}\}^L_{\ell=2},\bm{r}_{k+1}, \texttt{P\_hist}, \texttt{Re}, \bm{e})$
		\State 
		$\{\bm{A}_{\ell}\}^L_{\ell=2},\{\bm{P}_{\ell}\}^L_{\ell=1} \gets \texttt{PathCoverAMG\_setup}(\bm{A},\bm{e})$
		\State
		\texttt{Re} $\gets$ \texttt{Re} $+1$
		\State
		\texttt{P\_hist}$\{\texttt{Re}\} \gets \{\bm{P}_{\ell}\}^L_{\ell=1}$
		\State
		$\bm{x}_k\gets\bm{x}_k+\bm{e}$
		\EndIf
		\State
		$k \gets k+1$
		\EndWhile
		\EndProcedure
		
	\end{algorithmic}
\end{algorithm}

\begin{algorithm}[h!]
	\caption{PC-$\alpha$AMG (for $\bm{b}=\bm{0}$)}\label{alg:adpAMG_0}
	\begin{algorithmic}[1]
		\Procedure{[ $\bm{x}$ ] = \texttt{AdaptiveAMG}($\bm{A},\bm{b},\bm{x},\texttt{tol},\texttt{max\_iter}, \texttt{threshold}$)}{}
		\State
		$\{\bm{A}_{\ell}\}^L_{\ell=2},\{\bm{P}_{\ell}\}^L_{\ell=1}\gets$\texttt{MWM\_setup}($\bm{A}$)
		\State
		$k \gets 1$
		\Comment Initialization for iteration number
		\State
		$\bm{x}_k \gets \bm{x}$
		\State
		$\bm{r}_k \gets \bm{b}-\bm{A} \bm{x}_{k}$
		\While {$k<$ \texttt{max\_iter} and $\|\bm{r}_k\|>\texttt{tol}$} 
		\State
		$\bm{x}_{k+1}\gets$ \texttt{V\_cycle}$(\bm{A},\bm{b},\bm{x}_k,\{\bm{P}_{\ell}\}^L_{\ell=1},\{\bm{A}_{\ell}\}^L_{\ell=2})$
		\State 
		Make $\bm{x}_k$ orthogonal to $\bm{1}$
		\Comment If we are solving Graph Laplacian
		\State
		$\bm{r}_{k+1}\gets \bm{b}-\bm{A} \bm{x}_{k+1}$
		\If{$\frac{\|\bm{r}_{k+1}\|}{\|\bm{r}_{k}\|}>\texttt{threshold}$}
		\State
		$\bm{e}\gets\frac{\bm{x}_{k+1}-\bm{0}}{\|\bm{x}_{k+1}-\bm{0}\|}$
		\State 
		$\{\bm{A}_{\ell}\}^L_{\ell=2},\{\bm{P}_{\ell}\}^L_{\ell=1} \gets \texttt{PathCoverAMG\_setup}(\bm{A},\bm{e})$
		\EndIf
		\State
		$k \gets k+1$
		\EndWhile
		\EndProcedure
		
	\end{algorithmic}
\end{algorithm}

\section{Numerical Results}\label{sec:num}
In this section, we conduct numerical experiments on different types of graph Laplacians, some of them are related to discrete PDEs and the rest of them are from real-world networks.

We compare the performance of V-cycle UA-AMG using MWM coarsening scheme and Gauss-Seidel smoother with proposed \cref{alg:adpAMG_b} and \cref{alg:adpAMG_0}.  For \cref{alg:adpAMG_b}, we choose two values for \texttt{threshold}, which results in two different re-setup strategies. In the first case (denoted as \cref{alg:adpAMG_b}(1)), we choose $\texttt{threshold}=10^{-6}$, which essentially forces the multilevel hierarchy to be re-built after each iteration. Since the newly built hierarchy is specifically for eliminating the current smooth error, we expect that re-setup at each step would give the best performance in terms of iteration counts. However, re-setup every iteration is computationally expensive.  Therefore, we also consider another case where we set $\texttt{threshold}=0.4$ (denoted as \cref{alg:adpAMG_b}(2)) and, in this way, re-setup will not be triggered every iteration.  This choice balances the re-setup times and error reduction efficiency, which potentially reduces the total CPU time. 
For \cref{alg:adpAMG_0}, we simply use $\texttt{threshold}=0.5$ which seems to give the best performance in our numerical experiment.

For other parameters used in \cref{alg:adpAMG_b}, we set $\texttt{numW}=2$, $\texttt{iterW}=7$ and $\texttt{iterV}=4$ to generate the approximated smooth error. That is, we use W-cycle for the first two re-setups in order to get a good approximated error at the early stage and, later on, just use a few V-cycles to approximate the error. The relaxation scheme is Gauss-Seidel and the solving step adopts a basic multigrid V-cycle.


We chose to test our model problems with three different right-hand-sides: a low-frequency $\bm{b}\in \mathbb{R}^{|V|}=[1,1,\dots,1-|V|]^T$, 
a randomly generated $\bm{b}$ with zero sum, and $\bm{b}=\bm{0}$. For all the experiments below, the solver stops when the residual reaches the pre-described tolerance.  Basically, we use a relatively small tolerance $1e-8$ for the scalability tests in \cref{sec:regular-grid} and \cref{sec:ring-graph} and more practical choice $1e-6$ for real world graphs.
``--" in any columns of tables means the number of iterations exceeds $2500$. 

Numerical experiments are conducted using an $3$GHz Intel
Xeon `Sandy Bridge' CPU with $256$ GB of RAM.  The software used is an
algebraic multigrid package written by the authors and implemented in
Matlab. We report the following results in each table. The ``ConvR"
under MWM column is the algebraic average of convergence rate for the
last $10$ iterations, and since we restart whenever convergence rate at current iteration or on average is
above a threshold, we omit this factor for our adaptive algorithms. The
``Iter" column reports number of iterations needed for the
residual to reach certain tolerance, and ``Re'' records the number of re-setups needed specifically 
for \cref{alg:adpAMG_b} with high \texttt{threshold} and \cref{alg:adpAMG_0} (since \cref{alg:adpAMG_b}(1) re-setups every time, $\text{Re} = \text{Iter} -1$). ``OC'' is the operator
complexity defined in~\cite{Ruge.J;Stuben.K1987}, i.e. the ratio of
the total number of nonzeros of matrices on all levels divided by the
number of nonzeros of the finest level matrix. The operator complexity
in \cref{alg:adpAMG_b} is calculated as the algebraic average
of the operator complexities of all generated hierarchies. Finally, we report ``$t$" as the total
CPU times in seconds.

\subsection{Tests on Graphs Corresponding to Regular Grids} \label{sec:regular-grid}
We first tested the performance of the algorithms on unweighted graph Laplacians of 2D regular uniform grids in \cref{tab:grids_lfb}-\cref{tab:grids}. This is related to solving a Poisson equation with Neumann boundary condition on a 2D square.


\begin{table}[h]
	\begin{center}
		{\normalsize
			\caption{Performance of Solving Regular Grids with Low-Frequency $\bm{b}$, tol=1e-8}\label{tab:grids_lfb}
			\begin{tabular}{|c||c|c|c||c|c||c|c|c|} \hline \hline
				& \multicolumn{3}{|c||}{UA-AMG w/MWM} & \multicolumn{2}{|c||}{\cref{alg:adpAMG_b}(1)}& \multicolumn{3}{|c|}{\cref{alg:adpAMG_b}(2)} \\ \hline
				$n$ & Iter &  ConvR & OC& Iter  & OC & Iter &  Re  & OC \\ \hline \hline
				$2^{12}$ & 278 & 0.948  & 2.07 &  7 & 2.10 & 13 & 6  & 2.10 \\
				$2^{14}$ & 878 & 0.979  & 2.08 &  9  & 2.12 & 15 & 7 & 2.11	\\
				$2^{16}$ & 1960 & 0.990  & 2.09  & 10  &  2.12 & 19 & 9  & 2.13\\		
				$2^{18}$ &  --   & 0.996   & 2.09 &   11   & 2.14 & 21 & 10  & 2.15\\
				\hline \hline
			\end{tabular}
		}
	\end{center}
\end{table}

\begin{table}[h]
	\begin{center}
		{\normalsize
			\caption{Performance of Solving Regular Grids with Zero-Sum Random $\bm{b}$, tol=1e-8}\label{tab:grids_zsb}
			\begin{tabular}{|c||c|c|c||c|c||c|c|c|} \hline \hline
				& \multicolumn{3}{|c||}{UA-AMG w/MWM} & \multicolumn{2}{|c||}{\cref{alg:adpAMG_b}(1)}& \multicolumn{3}{|c|}{\cref{alg:adpAMG_b}(2)} \\ \hline
				$n$ & Iter &  ConvR  & OC& Iter  & OC & Iter &  Re  & OC \\ \hline \hline
				$2^{12}$ & 269 & 0.949  & 2.07 &  7  & 2.10 & 14 & 6  & 2.10	\\
				$2^{14}$ & 625 & 0.979  & 2.08 &  8  & 2.15 & 16 & 7  & 2.14	\\
				$2^{16}$ & 1383 & 0.991  & 2.09  & 10  &  2.16 & 18 & 8  & 2.15\\		
				$2^{18}$ &  --   & 0.996 & 2.10 &   11  &  2.17 & 20 & 10  & 2.16\\
				\hline \hline
			\end{tabular}
		}
	\end{center}
\end{table}

\begin{table}[h]
	\begin{center}
		{\normalsize
			\caption{Performance of Solving Regular Grids with $\bm{b}=\bm{0}$, tol=1e-8}\label{tab:grids}
			\begin{tabular}{|c||c|c|c||c|c|c|} \hline \hline
				& \multicolumn{3}{|c||}{UA-AMG w/MWM} & \multicolumn{3}{|c|}{\cref{alg:adpAMG_0}} \\ \hline
				$n$ & Iter &  ConvR & OC& Iter &  Re & OC \\ \hline \hline
				$2^{12}$ & 129 & 0.945  & 1.97 & 16 & 4 & 2.12\\
				$2^{14}$ & 194 & 0.964  & 1.99 &  16  & 4  & 2.13 	\\
				$2^{16}$ & 258 & 0.983 & 1.99  & 16 & 4  & 2.15\\		
				$2^{18}$ &  375   & 0.991   & 2.00 &   16  & 4  & 2.17\\
				\hline \hline
			\end{tabular}
		}
	\end{center}
\end{table}

\begin{figure}[h!]
	\begin{center}
		\begin{tabular}{cc}
			\includegraphics[width=0.45\textwidth]{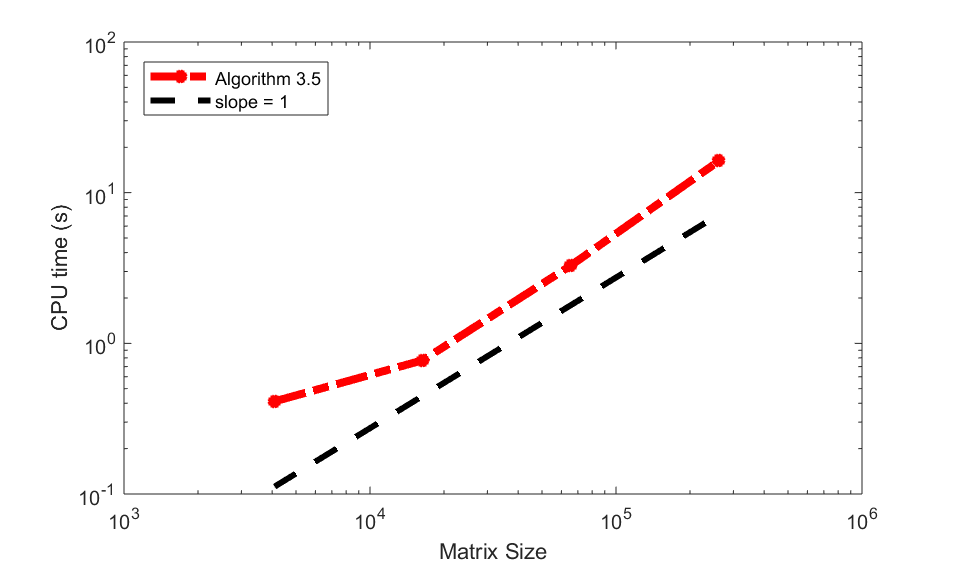}
			&  \includegraphics[width=0.45\textwidth]{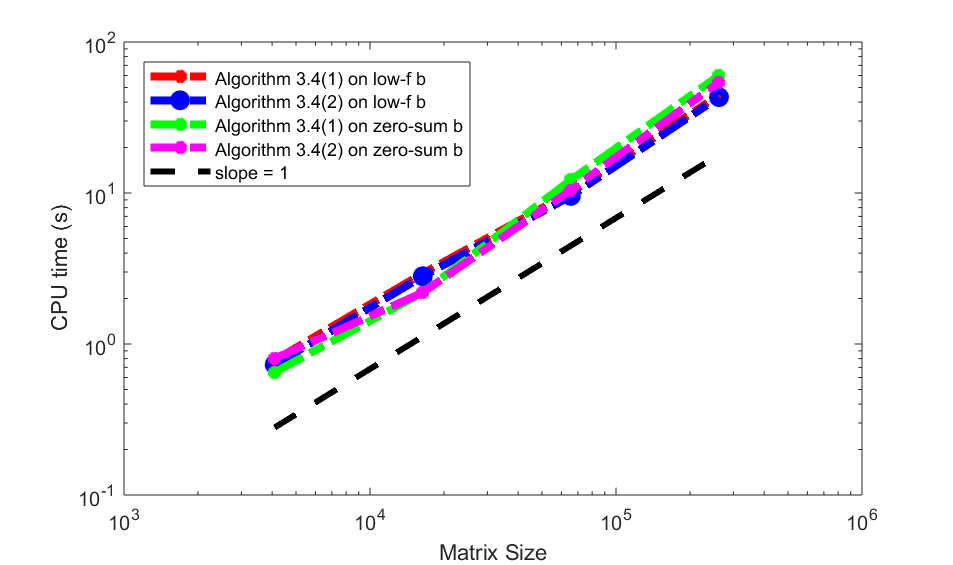}
		\end{tabular}
		\caption{CPU Time Elapsed for Solving Regular Grids with $\bm{b}=\bm{0}$ (left) and non-zero $\bm{b}$ (right)}\label{fig:grids_cpu}
	\end{center}
\end{figure}
Notice that for regular grids, while the numbers of iterations for regular V-cycle UA-AMG in all cases grow rapidly and quickly exceed 2500, number of iterations of \cref{alg:adpAMG_b} in \cref{tab:grids_lfb} and \cref{tab:grids_zsb} is nearly uniform and in \cref{tab:grids} is uniform. Total CPU time is plotted in \cref{fig:grids_cpu} and we can see that the total CPU time of PC-$\alpha$AMG increases nearly linearly with respect to the matrix size for both homogeneous $\bm{b}$ and non-zero $\bm{b}$ (the line with slope $1$ is added for reference). 
The growth rate of the CPU time is slightly slower than linear for small $n$, which is probably due to the fact that overhead cost is more pronounced than actually computing cost when $n$ is small.  But for large $n$, we can see the nearly linear growth asymptotically, which demonstrate that the computational cost of our PC-$\alpha$AMG is nearly optimal.  In addition, although compared to \cref{alg:adpAMG_b}(1), \cref{alg:adpAMG_b}(2) takes more iterations to converge, it needs slightly less CPU time because it re-setups less times than \cref{alg:adpAMG_b}(1).  This justifies our choices of different $\texttt{threshold}$.

\subsection{Tests for Ring Graphs} \label{sec:ring-graph}
\begin{figure}[h!]
	\begin{center}
			\includegraphics[width=0.6\textwidth]{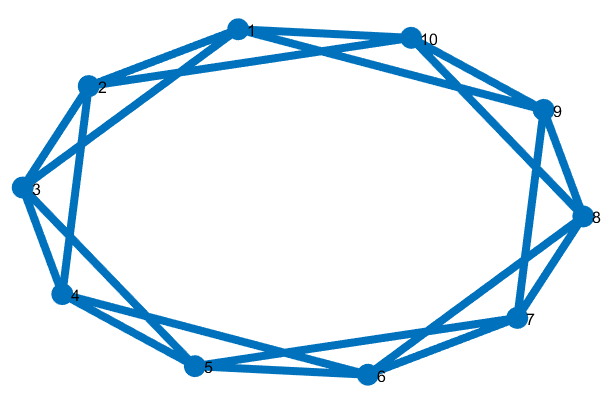}
		\caption{Example of a ring graph with $|V|=10$}\label{fig:ring}
	\end{center}
\end{figure}

For the second example, we use Watts Strogatz~\cite{watts1998collective} model and set the rewiring probability $ \beta = 0$ and the mean node degree to be 4 in order to produce ring graphs as in \cref{fig:ring}.  The condition numbers of the graph Laplacians of the ring graphs also grow rapidly when the size of the graphs increases.

\begin{table}[h]
	\begin{center}
		{
			\caption{Performance of Solving Ring Graphs with Low-Frequency $\bm{b}$, tol=1e-8}\label{tab:ring_lfb}
			\begin{tabular}{|c||c|c|c||c|c||c|c|c|} \hline \hline
				& \multicolumn{3}{|c||}{UA-AMG w/MWM} & \multicolumn{2}{|c||}{\cref{alg:adpAMG_b}(1)}& \multicolumn{3}{|c|}{\cref{alg:adpAMG_b}(2)} \\ \hline
				$n$ & Iter &  ConvR  & OC& Iter  & OC & Iter &  Re  & OC \\ \hline \hline
				$2^{10}$ & 149  & 0.894  & 1.59 & 7  & 1.57 & 14 & 6  & 1.57\\
				$2^{12}$ & 593  & 0.974  & 1.59 & 10  & 1.59 & 19 & 9  & 1.59\\
				$2^{14}$ & 2292 & 0.993  & 1.60 &  12   & 1.60 & 21 & 11 & 1.59	\\
				$2^{16}$ & -- & 0.998    & 1.60 & 18  & 1.60 & 24 & 15  & 1.60\\		
				\hline \hline
			\end{tabular}
		}
	\end{center}
\end{table}

\begin{table}[h]
	\begin{center}
		{
			\caption{Performance of Solving Ring Graphs with Zero-Sum Random $\bm{b}$, tol=1e-8}\label{tab:ring_zsb}
			\begin{tabular}{|c||c|c|c||c|c||c|c|c|} \hline \hline
				& \multicolumn{3}{|c||}{UA-AMG w/MWM} & \multicolumn{2}{|c||}{\cref{alg:adpAMG_b}(1)}& \multicolumn{3}{|c|}{\cref{alg:adpAMG_b}(2)} \\ \hline
				$n$ & Iter &  ConvR  & OC& Iter  & OC & Iter &  Re  & OC \\ \hline \hline
				$2^{10}$ & 144   & 0.10 & 1.56 & 8  & 1.57 & 14 & 6  & 1.58\\
				$2^{12}$ & 599  & 0.974   & 1.59 & 10  &  1.60 & 18 & 10  & 1.60\\
				$2^{14}$ & 2250 & 0.993  & 1.60 &  14   & 1.60 & 23 & 12  & 1.60	\\
				$2^{16}$ & -- & 0.998  & 1.60  & 21   & 1.60 & 25 & 18  & 1.60\\		
				\hline \hline
			\end{tabular}
		}
	\end{center}
\end{table}

\begin{table}[h]
	\begin{center}
		{\normalsize
			\caption{Performance of Solving Ring Graphs with $\bm{b}=\bm{0}$, tol=1e-8}\label{tab:ring}
			\begin{tabular}{|c||c|c|c||c|c|c|} \hline \hline
				& \multicolumn{3}{|c||}{UA-AMG w/MWM} & \multicolumn{3}{|c|}{\cref{alg:adpAMG_0}} \\ \hline
				$n$ & Iter &  ConvR  & OC& Iter &  Re  & OC \\ \hline \hline
				$2^{10}$ & 67  & 0.881  & 1.50 & 15 & 4  & 1.61\\
				$2^{12}$ & 311  & 0.980  & 1.50 & 15 & 4  & 1.62\\
				$2^{14}$ & 858 & 0.996   & 1.51 &  15  & 4   & 1.63 	\\
				$2^{16}$ & 995 & 0.997   & 1.51  & 15 & 4  & 1.63\\		
				\hline \hline
			\end{tabular}
		}
	\end{center}
\end{table}

\begin{figure}[h!]
	\begin{center}
		\begin{tabular}{cc}
			\includegraphics[width=0.45\textwidth]{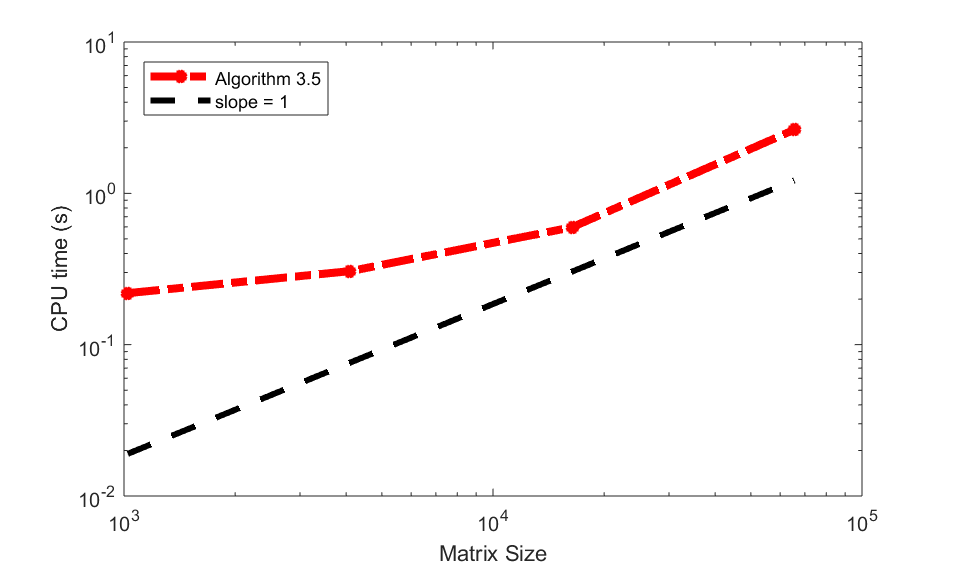}
			&  \includegraphics[width=0.45\textwidth]{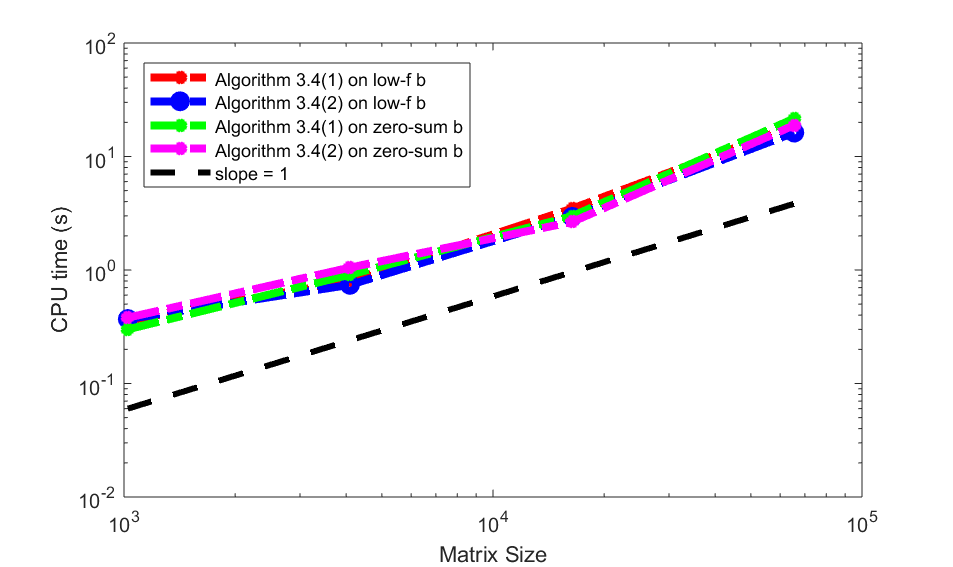}
		\end{tabular}
		\caption{CPU Time Elapsed for Solving Ring Graphs with $\bm{b}=\bm{0}$ (left) and non-zero $\bm{b}$ (right)}\label{fig:ring_cpu}
	\end{center}
\end{figure}


The results are present in \cref{tab:ring_lfb}, \cref{tab:ring_zsb} and \cref{tab:ring}. Like the performance of solving the regular grids, PC-$\alpha$AMG method on ring graphs requires
small numbers of iterations and re-setups to converge, while the standard V-cycle UA-AMG eventually cannot converge within $2500$ iterations for large graphs. 
In \cref{fig:ring_cpu}, we report the total CPU time for both zero and nonzero $\bm{b}$.  The results are similar to the regular grids' results.  When matrix size is small, the CPU time grows slower than linear.  Asymptotically, the CPU time grows nearly linear with respect to $n$, which demonstrate the nearly optimal computational complexity of our PC-$\alpha$AMG.  Moreover, compared to \cref{alg:adpAMG_b}(1), \cref{alg:adpAMG_b}(2) still takes slightly shorter CPU time for ring graphs.
\subsection{Tests for Real World Graphs}
Besides the graphs generated above, we also tested real world graphs from Stanford large network datasets collection~\cite{snapnets} and from the University of Florida Sparse Matrix Collection (UF)~\cite{davis2011university}. We selected graphs that are ill-conditioned and have relatively higher density.  Those graphs are quite challenging for standard AMG method.

We pre-processed the graphs as follows. The largest connected component of each graph is extracted, any self-loops from the extracted component are discarded, and edge weights of the component are modified to be their absolute values to satisfy the requirements of path cover finding algorithm. We also made the largest connected components undirected if they are originally directed. In \cref{tab:snapnets} and \cref{tab:UF}, the basic information of pre-processed graphs collected from two sources are presented. 

\begin{table}[htp]
	\caption{Largest connected components of the networks from Stanford large network datasets collection}\label{tab:snapnets}
	\begin{center}
	{\small
				\begin{tabular}{|c||c|c|l|} 
			\hline \hline
			&   n   &   nnz &   Description  \\ \hline \hline
			com-DBLP         & 3.17080e5    &  2.41681e6 &  DBLP collaboration network\\
			web-NotreDame &  3.25729e5  &  1.09011e6  &  Web graph of Notre Dame \\ 
			amazon0601     &  4.03364e5   &  5.28999e6  &  Amazon product co-purchasing network \\
			\hline \hline
		\end{tabular}
		}
	\end{center}
\end{table}%


\begin{table}[htp]
	\caption{Largest connected components of the networks from the University of Florida sparse matrix collection (UF)}\label{tab:UF}
	\begin{center}
		{\scriptsize
			\begin{tabular}{|c||c|c|l|} 
				\hline \hline
				&   n   &   nnz &   Description  \\ \hline \hline
				333SP & 3.71282e6 & 2.22173e7  & 2-dimensional FE
				triangular meshes	\\
				belgium\_osm &  1.44129e6 & 3.09994e6  & Belgium street network\\
				M6 &  3.50177e6   & 2.10038e7  & 2-dimensional FE
				triangular meshes\\
				NACA0015 &  1.03918e6   & 6.22963e6  & 2-dimensional FE
				triangular meshes\\
				netherlands\_osm &  2.21669e6   & 4.88247e6  & Netherlands street network\\
				packing-500x100x100-b050 &  2.14584e6   & 3.49765e7  & DIMACS Implementation Challenge\\
				roadNet-CA &  1.95703e6   & 5.52078e6  & California road network\\
				roadNet-PA &  1.08756e6   & 3.08303e6  & Philadelphia road network\\
				roadNet-TX &  1.35114e6   & 3.75840e6 & Texas road network\\
				fl2010	 &  4.84466e5   & 2.83072e6  & Florida census 2010  \\
				as-Skitter & 1.69642e6 & 2.21884e7 & Autonomous systems by Skitter\\		
				hollywood-2009 &  1.06913e6   & 1.13682e8  & Hollywood movie actor network\\
				\hline \hline
			\end{tabular}
		}
	\end{center}
\end{table}

\begin{table}[h]
	\begin{center}
		{
			\caption{Performance of Solving Graphs Collected from UF and Stanford with Low-Frequency $\bm{b}$, tol=1e-6}\label{tab:UFST_lbf}
			\begin{tabular}{|c||c|c|c||c|c||c|c|c|} \hline \hline
				& \multicolumn{3}{|c||}{UA-AMG w/MWM} & \multicolumn{2}{|c||}{\cref{alg:adpAMG_b}(1)}& \multicolumn{3}{|c|}{\cref{alg:adpAMG_b}(2)} \\ \hline
				 & Iter &  ConvR & OC& Iter & OC & Iter  & Re  & OC \\ \hline \hline
				 \multicolumn{9}{|c|}{UF large network datasets collection}\\
				 \hline \hline
				333SP & -- & 0.997  & 1.89 &  9   & 2.01 &  14 & 7 & 2.08	\\
				belgium\_osm &  1629   & 0.996   & 1.99 &   11   & 2.02 & 15 & 9  & 2.02\\
				M6 &  --   & 0.997   & 1.86 &  10  & 2.11 & 15 & 8  & 2.11\\
				NACA0015 &  --   & 0.995  & 1.86 &  9  & 2.10 & 14 & 7  & 2.10\\
				netherlands\_osm &  --   & 0.997  & 1.98 &   10  & 2.02 & 16 & 9  & 2.02\\
				packing &  --   & 0.999  & 1.06 &   11  & 2.46 & 19 & 10  & 2.46\\ 
				roadNet-CA &  878   & 0.991  & 2.05 &   8   & 2.08 & 14 & 7 & 2.08\\
				roadNet-PA &  1382   & 0.991   & 2.05 &   8  & 2.10 & 14 & 7  & 2.09\\
				roadNet-TX &  1424   & 0.994   & 2.04 &   9   & 2.08 & 14 & 7  & 2.08\\
				fl2010     &  --  & 0.998   & 1.83  
				& 9  & 2.19 & 15 & 7 & 2.19	\\
				as-Skitter &  -- & 0.998  & 1.21 & 10   & 3.13 & 19 & 8  & 3.14\\		
				hollywood-2009 &  --   & 0.999   & 1.01 & 5  & 3.17  &   11  & 3   & 3.18\\
				\hline \hline
				\multicolumn{9}{|c|}{Stanford large network datasets collection}\\
				\hline \hline
				com-DBLP       &  297  & 0.986  & 2.01
				& 4  &3.22 & 11 & 2  & 3.22 \\
				web-NotreDame   &  --   & 0.999   &   1.26 & 7 & 2.43  & 13 & 6  & 2.40\\
				amazon0601   &  -- & 0.998   & 1.58& 5  & 3.49 &  12 & 4   & 3.52 \\
				\hline \hline
			\end{tabular}
		}
	\end{center}
\end{table}

\begin{figure}[h!]
	\begin{center}
		\includegraphics[width=0.99\textwidth]{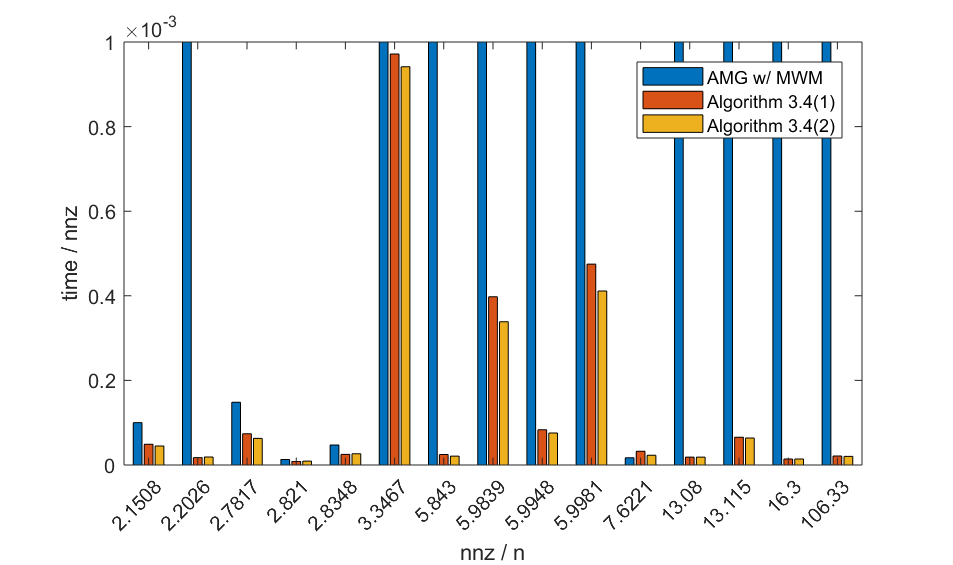}
		\caption{CPU Time Elapsed for Solving Real World Graphs with Low Frequency $\bm{b}$ using regular AMG and PC-$\alpha$AMG \cref{alg:adpAMG_b}}\label{fig:mat_lfb_cpu}
	\end{center}
\end{figure} 

\begin{table}[h]
	\begin{center}
		{
			\caption{Performance of Solving Graphs Collected from UF and Stanford with Zero-Sum Random $\bm{b}$, tol=1e-6}\label{tab:UFST_zsb}
			\begin{tabular}{|c||c|c|c||c|c||c|c|c|} \hline \hline
				& \multicolumn{3}{|c||}{UA-AMG w/MWM} & \multicolumn{2}{|c||}{\cref{alg:adpAMG_b}(1)}& \multicolumn{3}{|c|}{\cref{alg:adpAMG_b}(2)} \\ \hline
				& Iter &  ConvR & OC& Iter & OC & Iter  & Re  & OC \\ \hline \hline
				\multicolumn{9}{|c|}{UF large network datasets collection}\\
				\hline \hline
				333SP & -- & 0.997  & 1.89 &  9  & 2.09 &  6 & 1  & 2.08	\\
				belgium\_osm &   --   & 0.996  & 1.99 &   11   & 2.02 & 15 & 9  & 2.02\\
				M6 &  --   & 0.997  & 1.86 &  8  & 2.11 & 5 & 1  & 2.11\\
				NACA0015 &  1565   & 0.995   & 1.86 &  8  & 2.10 & 5 & 1  & 2.10\\
				netherlands\_osm &  --   & 0.997   & 1.98 &   12  & 2.02 & 17 & 11  & 2.02\\
				packing &  --   & 0.999  & 1.06 &   11  & 2.46 & 17 & 10  & 2.47\\ 
				roadNet-CA &  1308   & 0.994  & 2.08 &   8  & 2.08 & 15 & 7  & 2.08\\
				roadNet-PA &  970   & 0.991   & 2.05 &   8  & 2.09 & 14 & 6 & 2.08\\
				roadNet-TX &  1168   & 0.992   & 2.04 &   9   & 2.08 & 14 & 7  & 2.08\\
				fl2010     &  --  & 0.998  & 1.83  
				& 8  & 2.19 & 16 & 7  & 2.19	\\
				as-Skitter &  -- & 0.998  & 1.21 & 10   & 3.04 & 17 & 7  & 3.06\\		
				hollywood-2009 &  --   & 0.999  & 1.01 & 7  & 3.17  &   13  & 5   & 3.18\\
				\hline \hline
				\multicolumn{9}{|c|}{Stanford large network datasets collection}\\
				\hline \hline
				com-DBLP       &  573  & 0.987  & 2.01
				& 4  &3.23 & 11 & 3  & 3.22 \\
				web-NotreDame   &  --   & 0.999    &   1.26 & 7 & 2.47  & 15 & 6  & 2.56\\
				amazon0601   &  -- & 0.998  & 1.58& 6  & 3.49 &  10 & 4   & 3.50 \\
				
				\hline\hline
			\end{tabular}
		}
	\end{center}
\end{table}

\begin{figure}[h!]
	\begin{center}
		\includegraphics[width=0.99\textwidth]{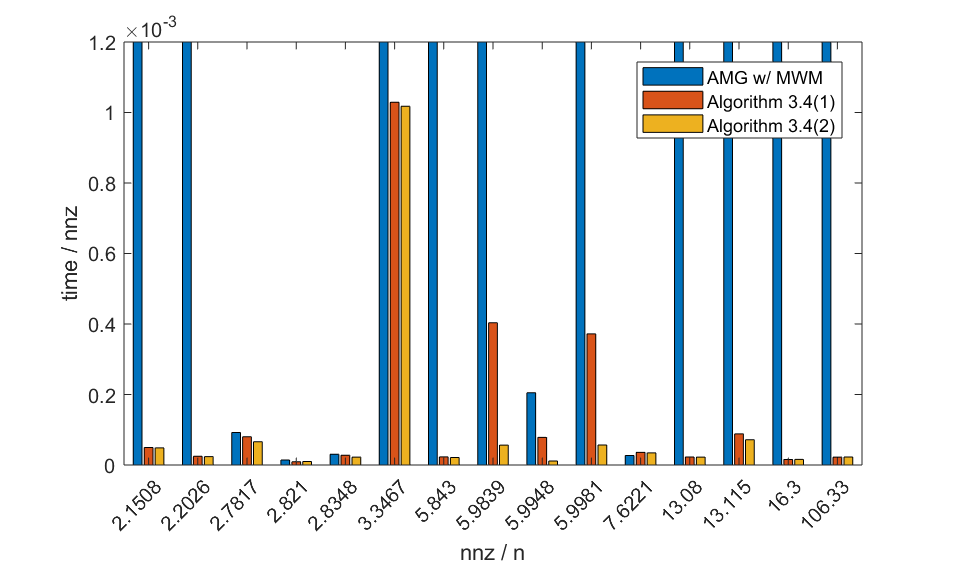}
		\caption{CPU Time Elapsed for Solving Real World Graphs with Randomly Generated Zero-Sum $\bm{b}$ using regular AMG and PC-$\alpha$AMG \cref{alg:adpAMG_b}}\label{fig:mat_zsb_cpu}
	\end{center}
\end{figure}

The results are presented in \cref{tab:UFST_lbf} and \cref{tab:UFST_zsb}.  As we can see, 
 for those graph Laplacians corresponding to more complicated (real-life)
graphs, whose properties are far from standard,  
standard V-cycle UA-AMG method struggles to converge and actually fail to converge for more than half of the graphs with low-frequency right-hand-side and zero-sum random right-hand-side. 
However, our PC-$\alpha$AMG converges within less than $20$ iterations for all tested graphs and building the coarse grid hierarchy need 10 re-setups on average, 
which demonstrates the effectiveness of PC-$\alpha$AMG.  Moreover, the average operator complexity in our adaptive algorithm is just slightly above 2, which suggests that our path covering aggregation scheme keeps the sparsity pattern on the coarse levels relatively well.  This fact, together with the usage of V-cycle,  makes our adaptive AMG method attractive for parallel computing and large-scale graphs.  


In \cref{fig:mat_lfb_cpu} and \cref{fig:mat_zsb_cpu}, we report the CPU time (more precisely, CPU time per number of non-zeros) for real world graphs.  We capped the height for all the cases where V-cycle UA-AMG did not converge within $2500$ steps. For low-frequency $\bm{b}$ (\cref{fig:mat_lfb_cpu}), we can observe that when the density (nnz/n) of the matrix is large, it is more likely that the V-cycle UA-AMG fails to converge within $2500$ iterations.  However, PC-$\alpha$AMG (both \cref{alg:adpAMG_b}(1) and \cref{alg:adpAMG_b}(2)) converges for all the cases and is faster than the regular AMG for all the tested graphs, especially for denser graphs.  Between \cref{alg:adpAMG_b}(1) and \cref{alg:adpAMG_b}(2), the CPU times are comparable while \cref{alg:adpAMG_b}(2) is slightly beter for some graphs.
For randomly generated zero-sum $\bm{b}$ (\cref{fig:mat_zsb_cpu}), the relationship between the convergence of V-cycle UA-AMG and the density of the matrices is more unpredictable. However, we can still observe that our PC-$\alpha$AMG outperforms V-cycle UA-AMG for all the tested graphs. 
 In addition, \cref{alg:adpAMG_b}(2) seems to be faster than \cref{alg:adpAMG_b}(1) for most of the tested graphs, which guides us to choose \cref{alg:adpAMG_b}(2) over \cref{alg:adpAMG_b}(1) for its flexibility.
 

Overall, our PC-$\alpha$AMG performs quite robust and efficient for graphs from real world applications, especially the highly ill-conditioned large-scale graphs.



\section{Conclusions}\label{sec:conc}
In this paper, we proposed PC-$\alpha$AMG, an adaptive UA-AMG algorithm, for solving linear systems based on optimal path cover.
The basic idea relies on adaptive construction of a multilevel
hierarchy as follows: (1) use a standard smoother to quickly reduce
the high frequency errors; (2) approximate the (algebraically) smooth
error on an adapted coarse grid using matching.  As the error changes
during the iterations, the second step may require a
re-setup (constructing a new multilevel hierarchy) to efficiently
eliminate the current smooth errors.  We approximate the level set
of the smooth error using a path cover and aggregate along the paths
(because the error is constant along paths following the level
set). We then build a coarse space based on such aggregation contains
a good approximation of the smooth error. The numerical tests on
different model problems show that, after each re-setup, the
dominating low frequency errors are quickly damped (with damping
factor $<0.2$ on average).  Thus, the proposed algorithm effectively
eliminates the algebraically smooth errors using several multilevel
hierarchies and, according to our numerical experiments, scales nearly optimally with respect to the size of the testing matrices even when applying standard V-cycle and unsmoothed aggregation scheme.

For $\bm{b}=\bm{0}$ case (\cref{tab:grids} and \cref{tab:ring}),
on each test problem, uniform convergence is observed. The work for
generating a new multilevel hierarchy is of order
$\mathcal{O}(|V|\log(|V|))$, as the path covering algorithm
runs only on the fine level, which costs $\mathcal{O}(|V|\log(|V|))$.
We also note that in the numerical tests, the number of re-setups needed is small
(4-5 on average) and is independent of the size or the type of the
model problem considered. Notice that when $\bm{b}=0$, the exact error is known and
these benchmark problems are just to show how the PC-$\alpha$AMG works.  In addition, in this case, PC-$\alpha$AMG can also be used as standalone setup phase for traditional adaptive AMG methods.

In the case of a non-zero right-hand-side, iteration count increases
slightly with the matrix size, since we can only use an approximation
of the smooth error in this case.  Total CPU time scales nearly linearly according to the numerical tests. Solving graph Laplacians corresponding to
real world graphs requires flexibility in choosing when to re-build a
new multilevel hierarchy. With such flexible choice, \cref{alg:adpAMG_b} requires less than $20$ iterations to achieve specified tolerance and the number of re-setups remains relatively small, which results in faster CPU time compared to standard V-cycle UA-AMG.
This shows the robustness and generality of our adaptive algorithm. 
We see such behavior on a wide range of matrices (graphs) tested (\cref{tab:UFST_lbf}
and \cref{tab:UFST_zsb}).

While the proposed algorithms clearly have the robustness to be
practical, we would also like to mention a few future research plans
which could improve such robustness.  
In our opinion, it is crucial to design
aggregation algorithms which approximate the errors well when
$\bm{b}\neq 0$.  As we mentioned earlier, a viable approach for this
is to use the adaptive aggregations based on a posteriori error
estimates as proposed in~\cite{Wenfang}. Another enhancement of our consideration is to involve more advanced aggregations/cycles/solvers
to approximate the smooth error. A combination of such approaches has
the potential to provide robust multilevel algorithms for solving
linear systems with graph Laplacians.


\bibliographystyle{siamplain} \bibliography{references}

\end{document}